\documentclass[a4paper]{amsart}
\usepackage{a4wide}
\usepackage{enumerate}
\usepackage{mathrsfs}
\usepackage{color}
\usepackage{amsfonts,amssymb,amscd,amsmath,amsthm,xypic}
\usepackage{txfonts}
\usepackage{accents}

\newif\iflineno 
\linenofalse
\iflineno \usepackage[pagewise]{lineno}\fi
\usepackage{ifpdf}

\ifpdf
    \usepackage[pdftex]{hyperref}
   \usepackage[pdftex]{graphicx}
   \newcommand{\myext}{pdf_tex}
\else
   \usepackage[dvips]{hyperref}
   \usepackage[dvips]{graphicx} 
   \newcommand{\myext}{ps_tex}
\fi

\newtheorem{Thm}[equation]{Theorem}

\theoremstyle{definition}

\newtheorem{Rem}[equation]{Remark}
\newtheorem{Fact}[equation]{Fact}

\theoremstyle{remark}

\numberwithin{equation}{section}

\newcommand{\DEFEQ}{\coloneqq}
\newcommand{\bL}{\mathbb{L}}
\newcommand{\bJ}{\mathbb{J}}
\newcommand{\BP}{\mathbf{P}}
\newcommand{\BQ}{\mathbf{Q}}
\newcommand{\forces}{\Vdash}

\newcommand{\prep}{\mathbb{R}}
\newcommand{\al}[1]{\ensuremath{{\aleph_{#1}}} }          % aleph_n
\newcommand{\om}[1]{\ensuremath{{\omega_{#1}}} }          % omega_n
\newcommand{\qemph}[1]{``\emph{#1}''}
\newcommand{\proofclaim}[2]{
\begin{equation}\label{#1}
  \parbox{0.8\textwidth}{#2}
\end{equation}}
\newcommand{\proofclaimnl}[1]{
\begin{equation*}
  \parbox{0.8\textwidth}{#1}
\end{equation*}}

\begin{document}
\iflineno \pagewiselinenumbers\fi

\date{2011-05-28}

\title[An overview of the proof in ``Borel Conjecture and dual Borel Conjecture'']
      {An overview of the proof in\\ ``Borel Conjecture and dual Borel Conjecture''}
\author{Martin Goldstern}
\address{Institut f\"ur Diskrete Mathematik und Geometrie\\
 Technische Universit\"at Wien\\
 Wiedner Hauptstra{\ss}e 8--10/104\\
 1040 Wien, Austria}
\email{martin.goldstern@tuwien.ac.at}
\urladdr{http://www.tuwien.ac.at/goldstern/}
\author{Jakob Kellner}
\address{Kurt G\"odel Research Center for Mathematical Logic\\
 Universit\"at Wien\\
 W\"ahringer Stra\ss e 25\\
 1090 Wien, Austria}
\email{kellner@fsmat.at}
\urladdr{http://www.logic.univie.ac.at/$\sim$kellner/}
\author{Saharon Shelah}
\address{Einstein Institute of Mathematics\\
Edmond J. Safra Campus, Givat Ram\\
The Hebrew University of Jerusalem\\
Jerusalem, 91904, Israel\\
and
Department of Mathematics\\
Rutgers University\\
New Brunswick, NJ 08854, USA}%{4}
\email{shelah@math.huji.ac.il}
\urladdr{http://shelah.logic.at/}
\author{Wolfgang Wohofsky}
\address{Institut f\"ur Diskrete Mathematik und Geometrie\\
 Technische Universit\"at Wien\\
 Wiedner Hauptstra{\ss}e 8--10/104\\
 1040 Wien, Austria}
\email{wolfgang.wohofsky@gmx.at}
\urladdr{http://www.wohofsky.eu/math/}
\thanks{
We gratefully acknowledge the following partial support: US National Science
Foundation Grant No. 0600940 (all authors); US-Israel Binational Science
Foundation grant 2006108 (third author); FWF Austrian Science Fund grant P21651
and EU FP7 Marie Curie grant PERG02-GA-2207-224747 (second and fourth author);
FWF grant P21968 (first and fourth author); \"OAW Doc fellowship (fourth author).}

\begin{abstract}
This note gives an informal overview of the proof in our paper ``Borel Conjecture and Dual Borel Conjecture''.
\end{abstract}

\maketitle

In this note, we give a rather informal overview (including two diagrams) of the
proof given in our paper ``Borel Conjecture and Dual Borel
Conjecture''~\url{http://arxiv.org/abs/1105.0823} (let us call it the ``main
paper''). This overview was originally a section in the main paper (the section
following the introduction), but the referee found it not so illuminating, so 
we removed it from the main paper.

Since we  think the introduction may be helpful to some readers
 (as opposed to the
referee, who found the diagrams ``mystifying'') we preserve it in this form.

The overview is supposed to complement the paper, and not to be read
independently. In particular,  we refer to the main paper for references.

The emphasis of this note is on giving the reader some vague understanding, at
the expense of correctness of the claims (we point out some of the most blatant
lies).

\section{The general setup}\label{mysec:generalsetup}

We assume CH in the ground model.  We use a $\sigma$-closed $\al2$-cc
preparatory forcing $\prep$, which adds a generic ``alternating iteration'' (as
defined below) $\bar \BP=(\BP_\alpha,\BQ_\alpha)_{\alpha<\om2}$.  Moreover,
$\prep$ forces that $\bar  \BP$ is ccc.  The forcing notion to get BC+dBC is the
composition  $\prep* \BP_{\om2}$.

We say that $\bar P$ is an ``alternating iteration'' if $\bar
P=(P_\alpha,Q_\alpha)_{\alpha<\om2}$ is a forcing iteration of length
$\om2$ satisfying the following:
\begin{itemize}
  \item At every even step $\alpha$, ($P_\alpha$ forces that)
    $Q_\alpha$ is an ``ultralaver forcing'' (described below).
  \item At every odd step $\alpha$, ($P_\alpha$ forces that) 
    $Q_\alpha$ is a ``Janus forcing'' (described below).
  \item However, instead of using either a Janus or an ultralaver
    forcing, we are at any step allowed just to ``do nothing'', i.e., 
    set $Q_\alpha=\{\emptyset\}$.
	\item At a limit step $\delta$, we take ``partial countable support'' limits.
    (This means more or less:
    $P_\delta$ is a subset 
    of the countable support limit of $(P_\alpha)_{\alpha<\delta}$
    and contains $\bigcup_{\alpha<\delta} P_\alpha$ and 
    has some other natural properties.)
\end{itemize}

\section{Ultralaver forcing}\label{mysec:ultralaver}

Let $\bar D=(D_s)_{s\in\omega^{<\omega}}$ be a system of ultrafilters.  The
``ultralaver forcing'' $\bL_{\bar{D}}$ consists of trees $p$ with the following
property: For every node $s\in p$ above the stem the set of immediate
successors of $s$ is in $D_s$. So this is a $\sigma$-centered variant of Laver
forcing. Of course this forcing adds a naturally defined generic real, called
ultralaver real.

We will basically need two properties of ultralaver forcing: The first one is preservation of positivity:
\proofclaim{eq:ulpos}{$\bL_{\bar{D}}$ preserves Lebesgue outer measure
positivity of ground model sets.\footnotemark}\footnotetext{This is a lie, and
moreover a stupid (i.e., useless) lie. It is a lie, since we only get something
like: for one random over a specific model, we can find a system $\bar{D}$ such that
$\bL_{\bar{D}}$ preserves randomness. It is a useless lie, since preservation of
positivity is not enough anyway: We need a stronger property that is preserved under
proper countable support iterations.}
The second one is killing of smz sets:
\proofclaimnl{For every uncountable set $X$ in the ground model,
$\bL_{\bar{D}}$ forces that $X$ is non-smz.}
Actually, we should formulate this claim in a stronger form.
Let us first quote a result of Pawlikowski, which is 
essential for the part of our proof that shows BC:
\begin{Thm}\label{thm:pawlikowski}
  $X\subseteq 2^\omega$ is smz iff $X+F$ is null for every closed null set $F$.
  \\
  Moreover, for every dense $G_\delta$ set $H$ we can \emph{construct} (in
  an absolute way) a closed null set $F$ such that for every $X
  \subseteq 2^\omega$ with $X+F$ null there is $t\in 2^\omega$ with
  $t+X\subseteq H$. 
\end{Thm}

So we can actually show the following:
\proofclaim{eq:ulsmz}{We can construct from the ultralaver real in an absolute
way a (code for a) closed null set $F$ such that 
$X+F$ is (outer Lebesgue measure) positive for every uncountable
ground model set $X$.}

It is an easy exercise to show that
 Theorem~\ref{thm:pawlikowski} implies  the following fact.

\begin{Fact} 
Assume that   $\bar P = (P_\alpha, Q_\alpha:\alpha <\omega_2)$ is 
an  iteration with direct limit $P_\om2$ satisfying the following:
\begin{itemize}
\item For cofinally many $\alpha<\omega_2$, $Q_\alpha$ makes every old uncountable set 
non-smz.  % (For example, $Q_\alpha$ = Laver forcing.) 
% \item For all $\alpha$, $Q_\alpha$ preserves outer measure. 
\item  $P_\om2$ and even all quotients $P_\om2/P_\alpha$
 preserves Lebesgue  outer measure positivity. 
\item $P_\om2$ preserves $\aleph_1$  and satisfies the $\al2$-cc.
\end{itemize}
Then $P_\om2$ forces BC. 
\end{Fact}

\begin{Rem} 
It is well-known that both Laver reals and random reals preserve positivity.
As Laver forcing makes every old uncountable set non-smz, we conclude that 
a countable support iteration of length $\omega_2$ of Laver reals, 
or alternatively, a countable support iteration alternating
Laver with random reals,   forces BC. The latter 
iteration also forces the failure of dBC, since the random reals increase
the covering number of the null ideal, and every set smaller than this 
cardinal is sm.
\end{Rem}

\section{A preparatory forcing for a single step}\label{mysec:prepsingle}

Let us first describe how to generically create a single forcing, e.g., an
ultralaver forcing.  

Let $Q$ be a forcing, $M$ a countable transitive model, $P\in M$ a
subforcing 
of~$Q$.  We say that   $P$ is an $M$-complete subforcing of~$Q$,
if every maximal antichain $A\in M$ of $P$ is also a maximal antichain
in~$Q$. In this case every $Q$-generic filter over $V$ induces a $P$-generic filter
over $M$.

Let $M^x$ be a (countable) model, and $\bar D^x\DEFEQ(D^x_s)_{s\in\omega^{<\omega}}$
a system of ultrafilters in $M^x$.  This defines the
ultralaver forcing $Q^x=\bL_{\bar{D}^x}$ in $M^x$.  Given
any system $\bar D$ of ultrafilters (in $V$) such that each $D_s$ extends
$D^x_s$, then we can show that $Q^x$ is an $M^x$-complete subforcing of
the ultralaver forcing $Q\DEFEQ\bL_{\bar{D}}$ in $V$. We describe this by\footnote{Note the linguistic asymmetry here: A symmetric and more verbose variant 
would say ``$(M^x,Q^x)$ canonically embeds into $(V,Q)$''.}
\qemph{$(M^x,Q^x)$ canonically embeds into $Q$}.

So every $Q$-generic filter $H$ over $V$ induces a $Q^x$-generic filter
over $M^x$ which we call $H^x$. A trivial but crucial observation is the
following: When we evaluate the ultralaver real 
for $Q$ in $V[H]$ then we get the
same real as when we evaluate it for $Q^x$ in $M^x[H^x]$.
Of course $Q^x$ is not a complete subforcing of $Q$, just an $M^x$-complete
one:
 While  $Q^x$ is
just countable, therefore
 equivalent to Cohen forcing 
 (from the point of view of $V$), the real added by the $Q^x$-generic 
 $H^x$ is an ultralaver real (over $M^x$ as well as over $V$), and therefore does
 not add a Cohen real over $V$. 

We can define a preparatory forcing $\prep_\bL$ (for a single ultralaver forcing)
consisting of pairs $x=(M^x,Q^x)$ as above ($M$ in some $H(\chi^*)$, say), 
and ordered as follows: $y$ is stronger than $x$ if $M^x\in
M^y$ and ($M^y$ thinks that) $(M^x,Q^x)$ canonically embeds into $Q^y$. It is
not hard to see that $\prep_\bL$ is $\sigma$-closed, and adds in the extension a
generic ultralaver forcing $\BQ$ such that each $x$ in the generic filter
embeds into $\BQ$.

Let $G$ be $\prep_\bL$-generic (over $V$). So in $V[G]$, we know that $Q^x$ is
an $M^x$-complete subforcing of $\BQ$ for all $x\in G$. Let $H$ be
$\BQ$-generic (over $V[G]$). Then $H$ induces a $Q^x$-generic filter
over $M^x$ (which we call $H^x$) for all  $x\in G$. Let us repeat the trivial
observation: As above, each $M^x[H^x]$ will see the ``real'' ultralaver real
(i.e., the one of $V[G][H]$).

Note that ``canonical embedding'' is a form of approximation: If $x$ is in
the generic filter, we do not know everything about $\BQ$, but we know that
$Q^x\subseteq \BQ$ and that the maximal antichains of $Q^x$ in $M^x$ will be
maximal antichains in $\BQ$ as well.

We should at this stage mention another simple concept that will be used
several times: Given $(M^x,Q^x)$ and (in $V$) some ultralaver forcing $Q$ such
that $(M^x,Q^x)$ embeds into $Q$, we can take some countable elementary
submodel $N$ of $H(\chi^*)$ containing $(M^x,Q^x)$ and $Q$, and Mostowski-collapse
$(N,Q)$ to $y=(M^y,Q^y)$. Then $y$ is in $\prep_\bL$ and stronger than $x$.

\section{Janus forcing}\label{mysec:janus}

With  ``Janus forcings'' we denote a family of forcing notions (as in the case
of ``ultralaver forcings''). Every Janus forcing $\bJ$ is a subset of $H(\aleph_1)$
and has a countable ``core'' $\nabla$ (which is the same for every Janus
forcing) and  some additional ``stuffing''.\footnote{Actually, the definition
of Janus forcing additionally depends on a real parameter.  In our application,
we will use ultralaver forcings as even stages $\alpha$, and use a Janus
forcing defined from the ultralaver real in the stage $\alpha+1$. The
following claims about Janus forcings only hold for this situation; in
particular the ground-model sets mentioned have to live in the model before the
ultralaver forcing.}
The forcing $\nabla$ will add a generic  real (``Janus real'') 
coding a null set $Z_\nabla$.
The forcing  $\nabla$  will not be a complete subforcing of $\bJ$, but we will
require that all maximal antichains involved in the name $Z_\nabla$ are also
maximal in 
$\bJ$, so that $\bJ$ also adds a generic null set~$Z_\nabla$. 

The crucial combinatorial content of Janus forcings heavily relies on previous
work by Bartoszy\'nski and Shelah. 

Analogously to the case of ultralaver forcing, let $\prep_\bJ$ consist of pairs
$x=(M^x,Q^x)$ such that $M^x$ is a countable model and $Q^x$  a Janus forcing in
$M^x$.  Given $x\in\prep_\bJ$ and a Janus forcing $Q$ in $V$, we say that $x$
canonically embeds into $Q$ if $Q^x$ is an $M^x$-complete subforcing of $Q$,
and we set $y\leq x$ in $\prep_\bJ$ if $M^y$ thinks that $x$ canonically embeds
into $Q^y$.  Again $\prep_\bJ$ is a $\sigma$-closed forcing and adds
a generic object $\BQ$ that is (forced to be) a Janus
forcing; and for every $x$ in the $\prep_\bJ$-generic filter, $x$ canonically
embeds into $\BQ$.

As in the case of ultralaver forcing, the Janus real $Z_\nabla$  is absolute.

As opposed to the case of ultralaver forcing, every Janus forcing $Q^x$ in any
model $M^x$ is itself a Janus forcing in $V$; so (taking the collapse
of an elementary
submodel as above)  we trivially get:
\proofclaim{eq:janusctbl}{For every $x\in\prep_\bJ$ there is a $y\leq x$
  such that in $M^y$,  $Q^x$ is a countable Janus forcing (so
  in particular equivalent to Cohen forcing).}

A crucial property of Janus forcing is that we can make it into random
forcing as well:
\proofclaim{eq:janusrandom}{For every $x\in\prep_\bJ$ there is a $y\leq x$
  such that in  $M^y$,  $Q^y$ is forcing equivalent to random
  forcing.}

So here we see an important property of the preparatory forcing $\prep_\bJ$,
which might seem a bit paradoxical at first: ``Densely'', $\BQ$ seems to be
Cohen as well as random forcing. This two-faced behavior gives Janus forcing
its name; one could also describe this behavior as ``faking'' (faking 
to be Cohen and faking to be random).

Janus forcing is the forcing notion that replaces the Cohen real 
in the dBC part of the proof.
 The crucial point is:
\proofclaimnl{A countable Janus forcing makes every uncountable ground model set
of reals non-sm.}
Well, that is actually not much of a point at all: As Carlson has shown, 
this is achieved by a Cohen real, and obviously a countable Janus forcing is
equivalent to a Cohen real. And Carlson even showed: When adding 
a Cohen real, this  uncountable ground model set remains non-sm 
even after forcing with another forcing notion, provided this forcing
notion has precaliber~$\al1$. 

So what we actually claim for Janus forcing is the more ``explicit'' version of
our trivial (after Carlson) claim above. First, let us recall
the (trivially modifies) definition of strongly meager (sm):
notation:
\proofclaim{def:sm.demorgan}{
  $X$ is \emph{not} sm, iff there a null set $Z$ (called ``witness'') 
  such that $X+Z = 2^\omega$.
%  Let $Z$ be a null set. We say that ``$Z$ witnesses that $X$ is \emph{not} sm''
%	if there is no real~$t$ with $(X+t)\cap Z =
%  \emptyset$, or equivalently, if $X+Z = 2^\omega$.\\
%(Clearly $X$ is not sm iff there is a witness $Z$.)
}

So what we really claim is the following:
\proofclaim{eq:janussm}{
The canonical null set $Z_\nabla$ added by a countable Janus forcing 
has the property that  $X+Z_\nabla = 2^\omega$ for all uncountable
ground model sets $X$, and moreover $X+Z_\nabla = 2^\omega$ 
is preserved by every  subsequent $\sigma$-centered forcing.
}
Of course $Z_\nabla$ is interpreted as a code for a null set, not a concrete subset of
the reals (otherwise $X+Z_\nabla= 2^\omega$ could not hold when we add new reals).

Let us again note that it is important that we can construct the null set
$Z_\nabla$ (rather: the code) in an absolute way from the Janus real and
get~\eqref{eq:janussm}.

\section{The preparatory forcing for the iteration}\label{mysec:prep}

\begin{figure}[tb]
\centering
\scalebox{0.7}{%\def\svgwidth{10cm} 
%scalebox: scales everything;
%svgwidth: scales graphics but leaves text standardsize
\input{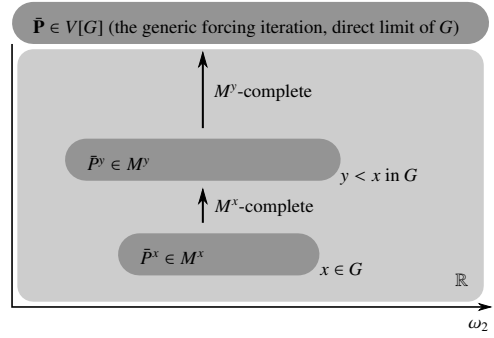}}
\caption{$G$ is the generic filter for the preparatory forcing notion $\prep$,
which adds the generic iteration $\bar\BP$.
The forcing that gives BC+dBC is $\prep*\BP_{\om2}$.
(Of course, in contrast to the impression given by the diagram, 
the set $M^x\cap\om2$, which is the ``domain'' of $\bar P^x$,
is not an interval.)}
\label{fig:prep}  
\end{figure}

The preparatory forcing $\prep$ that we will use will be similar to $\prep_\bL$
or to $\prep_\bJ$, but instead of ``approximating'' a single generic ultralaver
or Janus forcing, we approximate the alternating iteration $\bar\BP$ mentioned
in Section~\ref{mysec:generalsetup}.  So our preparatory forcing $\prep$ consists
of pairs $x=(M^x,\bar P^x)$, where $M^x$ is a countable model\footnote{Since we are
interested in iterations of length $\om2$, we cannot use transitive models,
that can only see ordinals $<\om1$.  Instead, we use ord-transitive models.}
(and subset of some fixed $H(\chi^*)$) and $\bar P^x$ is in $M^x$ an alternating
iteration.

Assume that $x\in \prep$ and that $\bar P$ is (in $V$) an alternating iteration.
As opposed to the case of a single ultralaver forcing, we now cannot formally
assume that $\bar P^x$ is a subset of $\bar P$, but there is a natural construction that
tries to give an $M^x$-complete embedding of $\bar P^x$ into $\bar P$. If this
construction works, we say that 
\qemph{$x$ is canonically embeddable into~$\bar P$},
and in that case we can treat $\bar P^x$ as subset of~$\bar P$.
So if $x$ is  canonically
embeddable into~$\bar P$,
 and $H$ is a $P_{\om2}$-generic filter over $V$ (which of course
induces $P_\alpha$-generic filters $H_\alpha$ for all $\alpha\leq\om2$), we can
get a canonical $P^x_\alpha$-generic filter over $M^x$ for every
$\alpha\in\om2\cap M^x$ (which we call $H^x_\alpha$).

We define the order in $\prep$ as above: For $x,y\in\prep$, we define
$y$ to be stronger than $x$, if $M^x\in M^y$ and $M^y$ thinks that $x$
canonically embeds into~$\bar P^y$.

Note that while $M^x$ thinks that $\bar P^x$ is an iteration of length $\om2$, in
$V$ (or in $M^y$ for $y\leq x$) the ``real domain'' of $ \bar P^x$ is just countable
(since it is a subset of $M^x$). 

As promised, one can show that $\prep$ is $\sigma$-closed and adds a generic
alternating iteration $\bar \BP$, and that $\bar \BP$ is ccc.  The final limit
$\BP_{\om2}$ is the direct limit of the $\BP_\alpha$
(and thus does not add any new   reals in the last stage). The  
intermediate stages satisfy CH, while $\BP_{\om2}$ forces $2^{\al0}=\al2$.  As might
be expected by now, each $x$ in the $\prep$-generic filter~$G$ canonically embeds
into~$\bar\BP$.  We will call the $\prep$-generic filter $G$. 
(The situation is illustrated in Figure~\ref{fig:prep}.) So if $H$ is
$\BP_{\om2}$-generic over $V[G]$, then we get canonical $P^x_{\om2}$-generic
filters $H^x$ for all $x\in G$; and the ``real'' ultralaver (and Janus) reals
calculated in $V[G][H]$ are the same as the ones ``locally'' calculated in
$M^x[H^x]$.

Given
any $x\in\prep$ we can construct (in $V$) an alternating
iteration $\bar P$ such that $x$ embeds into~$\bar P$
and such that $\bar P$ 
has either of the following two properties:
\begin{itemize}
  \item All Janus forcings are countable; at all stages $\alpha$ that
    are not in $M^x$ we ``do nothing''; and all
    limits $P_\delta$ are ``almost finite support over $x$''.
    (I.e., basically the limit is finite support, but we more or less add the 
    countably many elements of $P^x_\delta$.)
  \item All Janus forcings are equivalent to random forcing, and all
    limits $P_\delta$ are ``almost countable support over $x$''
    (basically we take all conditions in the countable support limit
    that are $x$-generic).
\end{itemize} 
The point is that these iterations behave more or less like finite (or
countable) support iterations; but we can still embed $x$ into them. 
For example, $M^x$ could
 think that $ \bar P^x$ is a countable support iteration, but we
 may still choose $\bar P$ to be an almost finite support iteration. 

``Behave more or less in the same way'' implies in particular in the first case
that any $P_\alpha$ is $\sigma$-centered: We iterate only countably many
forcings, since we do nothing outside $M^x$; the single forcings are
$\sigma$-centered (in the ultralaver case) and even countable in the Janus
case, and the (almost) finite support limits preserve $\sigma$-centeredness.

In the second case, we get preservation
of positivity (with respect to  outer Lebesgue
measure):
ultralaver as well as random forcings preserve positivity, and preservation is
preserved by (almost) countable support (proper) iterations.\footnote{Of
course, this is not true, rather we need an iterable property such as
preservation of random reals over models, etc. We do not get this stronger
property universally, we can just preserve a specific random; so
claim~\eqref{eq:janusB} is a lie, too.}

As above, we put the iteration $\bar P$ into a countable
elementary submodel; collapse it, and thus get:
\proofclaim{eq:janusA}{For all $x\in\prep$ there is a $y\leq x$
  such that ($M^y$ thinks that) $P^y_{\om2}$ is $\sigma$-centered and all
  Janus forcings are countable.}
\proofclaim{eq:janusB}{For all $x\in\prep$ there is a $y\leq x$
  such that ($M^y$ thinks that) $P^y_{\om2}$ preserves Lebesgue outer measure
  positivity.}

Let us again note that densely often we use finite support, but we also use 
countable support densely often.

\section{Why BC holds}\label{mysec:BC}

\newcommand{\mytextb}{Assume $\forces_\prep\forces_{\BP_\beta}X+F\subseteq Z$}
\newcommand{\mytextc}{$X^y\subseteq X$,}
\begin{figure}[tb]%
\centering%
\scalebox{0.8}{%\def\svgwidth{13cm} 
%scalebox: scales everything;
%svgwidth: scales graphics but leaves text standardsize
\input{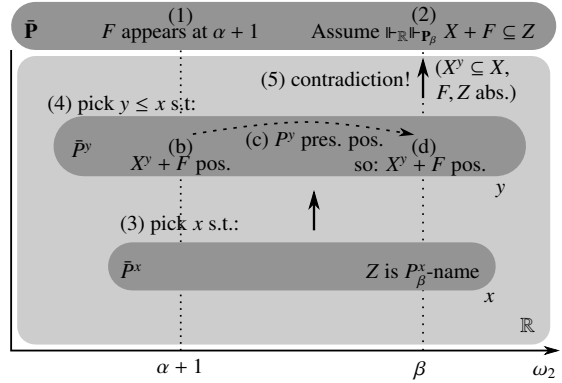}}%
\caption{The proof of BC.
}
\label{fig:BC_small}
\end{figure}
We want to show that BC is forced by $\prep*\BP_{\om2}$. Let $X$ be the name of
a set of reals of size~$\al1$.  Since $\BP_{\om2}$ has length $\om2$, we can
assume\footnote{Well, we can't. But we can do something similar, as will be
explained in the final section of the main paper.} that $X$ is in the ground
model $V$. We want to show BC, so we have to show  that $X$ is not smz. The
following is illustrated by Figure~\ref{fig:BC_small}.

\begin{enumerate}
  \item[(1)] Fix any ultralaver position $\alpha$. (Well, we fix $\alpha$ large
		enough to justify our assumption that $X\in V$.) We know that the
		ultralaver real that is added by $\BQ_\alpha$ (i.e., appears at stage
		$\alpha+1$) defines in an absolute way a (code for a) closed null set $F$.
  \item[(2)]
	 According to Theorem~\ref{thm:pawlikowski}, 
   it is enough to show that $X+F$ is non-null in the
   extension by $\prep*\BP_{\om2}$ (where the Borel code $F$ is evaluated in the extension).
   So assume towards a contradiction that
   $X+F$ is forced to be a subset of a null set (or rather, a Borel
   code) $Z$; this already has to
   happen\footnote{Each real (and in particular the Borel code $Z$) 
   in the $\BP_{\om2}$-extension already has to appear at some
	 stage $\beta<\om2$; and the statement ``$X+F\subseteq Z$''
	 is absolute.} at some 
   stage $\beta<\om2$. In other words: We assume (towards a contradiction)
   \[
    \forces_\prep \forces_{\BP_\beta}   X+F\subseteq Z. 
   \]
 \item[(3)] 
   Since $\BP_{\beta}$ is (forced to be) ccc, we
	 can find a very ``absolute'' (countable) name for $Z$; and we can find an $x\in\prep$
   that already calculates $Z$ correctly.\footnote{More formally: We find an
   $x\in\prep$ and a $P^x_\beta$-name $Z^x$ in $M^x$ such that $x$ forces (in
   $\prep$) that $\BP_\beta$ forces that $Z$ (evaluated by the
   $\BP_\beta$-generic) is the same as $Z^x$ (evaluated by the induced
   $P^x_\beta$-generic).}
 \item[(4)]
	 Now we construct (in $V$) a $y\leq x$ in $\prep$ (with $ \bar P^y_ \alpha$
	 proper)
	 that
	 satisfies~\eqref{eq:janusB}, and moreover such that 
	 $X^y\DEFEQ X\cap M^y\in M^y$
	 is uncountable in $M^y$ (we get this for free if $M^y$ is the collapse of an
	 elementary submodel $N$ with $X\in N$)
   In particular, ($M^y$ thinks that)
	 \begin{itemize}
	   \item[(a)] $P^y_\alpha$ is proper, thus preserves $\al1$,
       thus forces that $X^y$ is uncountable.
     \item[(b)] Therefore
		   $Q^y_\alpha$ forces 
		   that $X^y+F$ is positive (according to~\eqref{eq:ulsmz}).
	   \item[(c)] $P^y_\om2$ preserves positivity.
	   \item[(d)] Therefore $P^y_\beta$ forces that $X^y+F$ is positive\footnote{Here we even get positivity of $X^y+F$ where $F$ is evaluated in the intermediate extension of stage $\alpha+1$. However, we get the contradiction even if we just assume that $X^y+F$ is positive where $F$ is evaluated in the $P^y_\beta$-extension.}
       (and in particular not  a subset of $Z$).
	 \end{itemize}
  \item[(5)]
    This leads to the obvious contradiction:
		Let $G$ be $\prep$-generic over $V$ and contain $y$, and let $H_\beta$
		be $\BP_{\beta}$-generic over $V[G]$.  Then $H^y_\beta$ is
		$P^y_\beta$-generic over $M^y$, and therefore $M^y[H^y_\beta]$ thinks 
		that  some $x+f\in X^y+F$ is not in $Z$.
		But $X^y=X\cap M^y\subseteq X$, and the codes for $F$ and for $Z$
		are absolute (for $F$ since it is constructed in a canonical way
		from the ultralaver real, for $Z$ because we took care of it in step (3)).
\end{enumerate}

\section{Why dBC holds}\label{mysec:dBC}

The proof of dBC is similar, using $\sigma$-centeredness 
instead of positivity preserving, and a countable Janus forcing instead
of ultralaver forcing.

Let $X$ be the name of a set of reals of size~$\al1$.  Again, without loss of
generality\footnote{And again, this is a lie.} $X\in V$.  We want to show dBC,
so we have to show that $X$ is not sm.

\begin{enumerate}
  \item[(1)] Fix any Janus position $\alpha$ (large
		enough to justify our assumption that $X\in V$). We know that the
		Janus real that is added by $\BQ_\alpha$ (i.e., appears at stage
		$\alpha+1$) defines in an absolute way a (code for a) null set $Z_\nabla$.
  \item[(2)]
	 According to~\eqref{def:sm.demorgan}, it is enough to show that
   $Z_\nabla+X=2^\omega$
   in the extension by   $\prep*\BP_{\om2}$.
   So assume towards a contradiction that $Z_\nabla+X\neq 2^\omega$.
   This already happens at some stage $\beta<\omega$, i.e., we assume
   \[
    \forces_\prep\ \forces_{\BP_\beta}   r\notin Z_\nabla+X .
   \]
 \item[(3)] 
   Again, find $x$ such that $r$ is an ``absolute'' $P^x_\beta$-name.
 \item[(4)]
	 Now we construct (in $V$) a $y\leq x$ in $\prep$ that
	 satisfies~\eqref{eq:janusA}, and such that 
	 $X^y\DEFEQ X\cap M^y\in M^y$
	 is uncountable in $M^y$.
   In particular, ($M^y$ thinks that)
	 \begin{itemize}
	   \item[(a)] $P^y_\alpha$ is proper, thus preserves $\al1$, and so
       forces that $X^y$ is uncountable.
     \item[(b)]
       $Q^y_\alpha$ is (forced to be) a countable Janus forcing notion.
	   \item[(c)] $P^y_\om2$ is $\sigma$-centered.
	   \item[(d)] Therefore~\eqref{eq:janussm} implies
       that $P^y_\beta$ forces that $Z_\nabla+X^y=2^\omega$, in particular
       that $r\in Z_\nabla+X^y$.
	 \end{itemize}
  \item[(5)]
    As before, this leads to a contradiction.
\end{enumerate}
\end{document}